\RequirePackage{ifpdf}
\documentclass[12pt]{article}
\usepackage{amsthm}
\usepackage{amsmath}
\usepackage{amssymb}
\usepackage{subeqnarray}
\usepackage{cases}
\usepackage{latexsym}
\usepackage{graphicx}
\usepackage{fancybox,color}
\usepackage{longtable}
\usepackage{color}
\usepackage{lscape}
\usepackage{graphicx}
\usepackage{float}
\definecolor{lred}{rgb}{1,0.8,0.8}
\definecolor{lblue}{rgb}{0.8,0.8,1}
\definecolor{dred}{rgb}{0.6,0,0}
\definecolor{dblue}{rgb}{0,0,0.5}
\usepackage{url}

\newcommand{\blue}[1]{\begin{color}{blue}#1\end{color}}

\newcommand{\red}[1]{\begin{color}{red}#1\end{color}}
\newcommand{\green}[1]{\begin{color}{green}#1\end{color}}

\def\cA{{\cal A}}

\def\cM{{\cal M}}
\def\cP{{\cal P}}
 \def\R{{\cal R}}
\def\cS{{\cal S}}
\def\cT{{\cal T}}

\def\cS{{\cal S}} \def\S{{\cal S}}

\def\hY{\widehat{Y}}

\newcommand{\inprod}[2]{\langle #1 , #2 \rangle}

\def\norm#1{\|#1\|}
\def\abs#1{\left| #1\right|}
\def\cN{{\cal S}^n_{\geq 0}}
\def\cC{{\cal C}}
\def\cM{{\cal M}}

\def\cK{{\cal K}}

\def\hy{\hat{y}}
\def\hZ{\widehat{Z}}

\def\inprod#1#2{\langle#1, \, #2\rangle}

\def\diag{{\rm diag}}

\def\sig{\sigma}

\def\ty{\tilde{y}}
\def\tC{\widetilde{C}}

\def\tZ{\widetilde{Z}}
\def\tS{\widetilde{S}}
\def\tX{\widetilde{X}}

\def\cX{{\cal X}}

\def\hY{\widehat{Y}} \def\hy{\hat{y}}
\def\hZ{\widehat{Z}}
\def\hS{\widehat{S}}

\newtheorem{remark}{Remark}[section]

\newtheorem{thm}{Theorem}[section]  
\newtheorem{lem}[thm]{Lemma}        

\newtheorem{assump}{Assumption}

\newcommand{\T}{\mbox{\textrm{\tiny{T}}}}
\newcommand{\+}{\mbox{\tiny{$+$}}}

\def\mc{\multicolumn}

\def\inprod#1#2{\langle #1,\,#2\rangle}
\def\norm#1{\|#1\|}

\def\PP{\mbox{$D$}}
\def\DD{\mbox{$P$}}

\def\mc{\multicolumn}
\DeclareMathOperator*{\argmini}{arg \, min}

\newcommand{\argmin}{\argmini}

\newcommand{\clear}[1]{}
\begin{document}

\title{\bf SDPNAL$+$: A Majorized Semismooth Newton-CG Augmented Lagrangian Method for Semidefinite Programming with Nonnegative Constraints}

\author{
Liuqin Yang\thanks{Department of Mathematics, National University of Singapore, 10 Lower Kent Ridge Road, Singapore ({\tt yangliuqin@nus.edu.sg}).}, \,
Defeng Sun\thanks{Department  of  Mathematics  and  Risk  Management  Institute, National University of Singapore, 10 Lower Kent Ridge Road, Singapore ({\tt matsundf@nus.edu.sg}). This author's research was
supported in part by the Academic Research Fund under Grant R-146-000-149-112.}, \,
Kim-Chuan Toh\thanks{Department of Mathematics, National University of Singapore, 10 Lower Kent Ridge Road, Singapore
({\tt mattohkc@nus.edu.sg}). Research supported in part by 
the  Ministry of Education, Singapore, 
Academic Research Fund
under Grant R-146-000-194-112.}
}
\date{May 28, 2014}

\maketitle

\begin{abstract}
In this paper, we present a majorized semismooth Newton-CG  augmented Lagrangian method, called SDPNAL$+$, for semidefinite programming (SDP) with partial or full nonnegative constraints on the matrix variable. SDPNAL$+$ is a much enhanced version of SDPNAL introduced by Zhao, Sun and Toh [SIAM Journal on Optimization, 20 (2010), pp.~1737--1765] for solving generic SDPs. SDPNAL works very efficiently for nondegenerate SDPs but may encounter numerical difficulty for degenerate ones. Here we tackle this numerical difficulty by employing a majorized semismooth Newton-CG augmented Lagrangian method coupled with a convergent 3-block alternating direction method of multipliers  introduced recently by Sun, Toh and Yang [arXiv preprint arXiv:1404.5378, (2014)]. Numerical results for various large scale SDPs with or without nonnegative constraints show that the proposed method is not only fast but also robust in obtaining accurate solutions. It outperforms, by a significant margin, two other competitive publicly available  first order methods based   codes: (1) an alternating direction method of multipliers  based solver called SDPAD by Wen, Goldfarb and Yin [Mathematical Programming Computation, 2 (2010), pp.~203--230] and (2) a two-easy-block-decomposition hybrid proximal extragradient method called 2EBD-HPE by Monteiro, Ortiz and Svaiter [Mathematical Programming Computation,  (2013), pp.~1--48]. In contrast to these two codes, we are able to solve all the 95 difficult SDP problems arising from the relaxations of
quadratic assignment problems tested in SDPNAL to an accuracy of $10^{-6}$ efficiently, while SDPAD and 2EBD-HPE successfully solve 30 and 16 problems, respectively.
In addition, SDPNAL$+$ appears to be the only viable method currently available to solve large scale SDPs arising from rank-1 tensor approximation problems constructed by Nie and Wang [arXiv preprint arXiv:1308.6562, (2013)]. The largest rank-1 tensor approximation problem we solved (in about 14.5 hours) is {\tt nonsym(21,4)},  in which its resulting SDP problem has matrix dimension $n = 9,261$ and the number of equality constraints $m =12,326,390$. 
\end{abstract}



\section{Introduction}

Let $\cK$ be a pointed closed convex cone whose interior $int(\cK) \neq \emptyset$ and $\cP$ be a polyhedral convex cone in
a finite-dimensional Euclidean space $\cX$ such that $\cK \cap \cP$ is non-empty. For any cone $\cC \subseteq \cX$, we denote the dual cone of $\cC$ by $\cC^*$. For any closed convex set $\cC \subseteq \cX$, we denote the metric projection of $\cX$ onto $\cC$ by $\Pi_\cC(\cdot)$ and the tangent cone of $\cC$ at $X \in \cC$ by $\cT_{\cC}(X)$, respectively. We will make extensive use of the Moreau decomposition theorem in \cite{moreau1962decomposition}, which states that $X = \Pi_\cC(X) - \Pi_{\cC^*}(-X)$ for any $X\in \cX$ and any closed convex cone $\cC \subseteq \cX$. Let $\cS^n$ be the space of $n \times n$ real symmetric matrices and $\S^n_+$ be the cone of  positive semidefinite matrices in $\S^n$. In this paper, we focus on the case where $\cX = \cS^n$, $\cK = \cK^* = \S^n_+$. We are particularly interested in the case where $\cP = \cN$, the cone of  $n\times n$ real symmetric matrices whose elements are all nonnegative, though the algorithm which we will design later is also applicable to other cases. For any matrix $X \in \cS^n$, we use $X \succ 0$ to indicate that $X$ is a real symmetric positive definite matrix.

Consider the semidefinite programming (SDP) with an additional polyhedral cone constraint, which we name as SDP$+$:
\begin{eqnarray}
 ({\rm P}) \quad  \max \Big\{ \inprod{-C}{X} \mid \cA(X) = b,\; X \in \cK, \; X\in \cP \Big\},
  \label{eq-primal}
\end{eqnarray}
where $b\in\Re^m$ and $C \in \cX$ are given data, $\cA:\cX\rightarrow \Re^m$ is a given linear map whose adjoint is denoted as $\cA^*$.
Note that $\cP = \cX$ is allowed in \eqref{eq-primal}, in which case there is no additional polyhedral cone constraint imposed on $X$.
We assume that the matrix $\cA\cA^*$ is invertible, i.e., $\cA$ is surjective. The dual of (P) is given by
\begin{eqnarray}
 ({\rm D})\quad  \min \Big\{ \inprod{-b}{y} \mid \cA^*(y) + S + Z= C,\;
   S\in \cK^*, \; Z\in\cP^* \Big\}.
  \label{eq-dual}
\end{eqnarray}
The optimality conditions (KKT conditions) for (P) and (D) can be written as follows:
\begin{eqnarray} \label{eq-KKT}
\left\{ \begin{array}{l l}
  \cA(X) -b =  0, \quad   \cA^*(y) + S + Z - C = 0, &    \\[5pt]
   \inprod{X}{S}  =  0, \; X\in \cK,\; S\in \cK^*, \;
   \inprod{X}{Z}  =  0, \; X\in \cP,\; Z\in \cP^* .  &
 \end{array}
 \right.
\end{eqnarray}

In order for the KKT conditions (\ref{eq-KKT}) to have solutions, throughout this paper we make the following blanket assumption.
\begin{assump} \label{assump-feasible}
{\rm (a)} For problem {\rm (P)}, there exists a feasible solution $X_0 \in \cS^n_{+}$  such that
\begin{eqnarray}
\mathcal{A}(X_0) = b, \; X_0 \succ \mathbf{0}, \; X_0 \in \cP.
\end{eqnarray}
{\rm (b)} For problem {\rm (D)}, there exists a feasible solution $(y_0,S_0,Z_0) \in \Re^m \times \cS^n_{+} \times \cS^n$ such that
\begin{eqnarray}
\cA^*(y_0) + S_0 + Z_0 = C, \; S_0 \succ 0, \; Z_0 \in \cP^*.
\end{eqnarray}
\end{assump}
It is known from convex analysis (e.g, \cite[Corollary 5.3.6]{borwein2006convex}) that under Assumption \ref{assump-feasible}, the strong duality for (P) and (D) holds and the KKT conditions (\ref{eq-KKT}) have solutions.

For a given $\sig>0$, define the augmented Lagrangian function for the dual problem (D) as
follows:
\begin{eqnarray}
 L_{\sig}(y,S,Z; X) &=& \inprod{-b}{y} + \inprod{X}{\cA^*y + S + Z-C}
 + \frac{\sigma}{2} \norm{\cA^*y +S + Z-C}^2
  \nonumber \\[5pt]
 &=& \inprod{-b}{y} + \frac{\sigma}{2}
 \norm{\cA^*y + S + Z + \sigma^{-1} X - C}^2 - \frac{1}{2\sigma}
 \norm{X}^2,
\label{eq-dual-AL}
\end{eqnarray}
where $X \in \cX, y\in\Re^m, S \in \cK^*, Z \in \cP^*$.
We can consider the following inexact augmented Lagrangian method to solve (D). Specifically, given $\sigma_0 > 0$, $(y^0,S^0,Z^0) \in \Re^m \times \cK^* \times \cP^*$, perform the following steps at the $(k+1)$-th iteration:
\begin{subequations}
\begin{numcases}{}
   (y^{k+1},S^{k+1},Z^{k+1}) \approx \argmin
     \{ L_{\sig_k}(y,S,Z; X^k)\,\mid\, y\in\Re^m, \; S \in \cK^*, \; Z\in \cP^* \}, \label{eq-al-a}
   \\[5pt]
    X^{k+1} = X^k +  \sig_k (\cA^*y^{k+1}+S^{k+1}  + Z^{k+1}-C), \label{eq-al-b}
\end{numcases}
\end{subequations}
where $\sigma_k \in (0,+\infty)$, $k = 0,1,\ldots$.
For a general discussion on the convergence of the  augmented Lagrangian method for solving convex optimization problems and beyond, see \cite{rockafellar1976augmented,rockafellar1976monotone}.

Note that problem (P) can be reformulated as a standard SDP in the primal form by replacing the constraint $X \in \cP$ with two constraints $X-Y = 0$ and $ Y \in \cP$. In \cite{zhao2010newton},
SDPNAL introduced by Zhao, Sun and Toh is applied to solve such a reformulated problem. It works quite well for nondegenerate SDPs, especially those without the constraint $X \in \cP$. However, many of the tested SDPs (with the constraint $X \in \cP$) in \cite{zhao2010newton}
are degenerate and SDPNAL is  unable to solve those problems efficiently. Motivated by our
desire to overcome the aforementioned difficulty in solving degenerate SDPs and to improve the performance of SDPNAL, we present here a majorized semismooth Newton-CG augmented Lagrangian method by directly working on (P) instead of its reformulated problem. We call this new method SDPNAL$+$ since it is a much enhanced version of SDPNAL and it is designed for SDP$+$ problems (P).

The remaining parts of this paper are organized as follows. In Section \ref{sec-sub}, we introduce a majorized semismooth Newton-CG method for solving the inner minimization problems of the augmented Lagrangian method and analyze the convergence for solving these inner problems. Section \ref{sec-SDPNALP} presents the SDPNAL$+$ dual approach. Section \ref{sec-numerical} is on
numerical issues.  There we report numerical results for a variety of SDP$+$ and SDP problems. We make an extensive numerical  comparison with two other competitive  first order methods based codes:  (1) an alternating direction method of multiplier (ADMM)  based solver called SDPAD by Wen et al. \cite{wen2010alternating} and (2) a two-easy-block-decomposition hybrid proximal extragradient method called 2EBD-HPE by Monteiro et al. \cite{monteiro2013first}.   Numerical results   show that SDPNAL$+$ is both  fast and  robust in achieving accurate solutions.

For the first time, we are able to solve all the $95$ difficult SDP problems arising from
the relaxations of quadratic assignment problems (QAPs) tested in SDPNAL to an accuracy of $10^{-6}$ efficiently, while SDPAD and 2EBD-HPE successfully solve 30 and 16 problems, respectively. In addition, SDPNAL$+$ appears to be the only viable method currently available to solve large scale SDPs arising from rank-1 tensor approximation problems constructed by Nie and Wang \cite{nie2012regularization}. The largest rank-1 tensor approximation problem solved is {\tt nonsym(21,4)},  in which its resulting SDP problem has matrix dimension $n = 9,261$ and the number of equality constraints $m =12,326,390$. Finally, in order to demonstrate the power of the proposed majorized semismooth Newton-CG procedure, we list the numerical results by only running the  convergent ADMM with $3$-block constraints (ADMM+ in short) introduced by Sun et al. \cite{sun2014convergent}. As one may observe, although ADMM+ outperforms both SDPAD and  2EBD-HPE, it can still encounter numerical difficulty in solving some hard problems such as those arising from QAPs to high accuracy. The superior numerical performance of
SDPNAL$+$  over solvers based purely on first order methods such as
SDPAD and 2EBD-HPE clearly shows the necessity of exploiting   second order methods such as the semismooth Newton-CG method in order to solve
hard SDP$+$ and SDP problems to high accuracy efficiently. While
there has been a recent focus on using first order methods such as those
based on ADMM or accelerated proximal gradient methods to solve
structured convex optimization problems arising from machine learning
and statistics, the extensive numerical results we obtained here
for matrix conic programming
problems serve to demonstrate that  second order methods
with good local convergence property are essential, if used wisely, for mitigating
the inherent slow local convergence  of first order methods,
especially on difficult problems.



\section{A Majorized Semismooth Newton-CG Method for Inner Problems} \label{sec-sub}

Let $\sigma>0$ and $\tX \in \cS^n$ be fixed. In this section we will present a majorized semismooth Newton-CG method   for solving the following inner problems involved in the augmented Lagrangian method (\ref{eq-al-a}):
\begin{eqnarray}
   \min \Big\{ \phi(y,S,Z) := L_{\sig}(y,S,Z; \tX)\,\mid\, y\in\Re^m, \; S \in \cK^*, \; Z \in \cP^* \Big\}.
  \label{eq-al-sub-org}
\end{eqnarray}
Note that problem (\ref{eq-al-sub-org}) is the dual of the following problem:
\begin{equation}\label{suborgdual}
  \max \Big\{ \langle \,- C, X \rangle - \frac{1}{2
  \sigma}\|X - \tX \|^2 \,\mid \, \mathcal{A} (X) = b, \; X \in \cK, \; X \in \cP \Big\}.
\end{equation}
Since the objective function in (\ref{suborgdual}) is strongly concave, (\ref{suborgdual}) has a unique optimal solution. In order for its dual problem (\ref{eq-al-sub-org}) to have a bounded solution set, we need the following generalized Slater condition.

\begin{assump} \label{assump-slater}
There exists a positive definite matrix $X_0 \in \cS^n_{+} \cap relint(\cP)$ such that
\begin{eqnarray} \label{eq-slater}
\cA(\cT_{\cP}(X_0)) = \Re^m, \; X_0 \succ 0,
\end{eqnarray}
where $relint(\cP)$ denotes the relative interior of $\cP$.
\end{assump}

Note that in (\ref{eq-slater}), $\cT_{\cP}(X_0)$ is actually a linear subspace of $\cS^n$ as $X_0$ is assumed to be in the relative interior part of the polyhedral cone $\cP$.  When $\cP = \cS^n$, Assumption \ref{assump-slater} is equivalent to saying that
\begin{eqnarray}
\left\{ \begin{array}{l} \mathcal{A}: \mathcal{S}^n \rightarrow
\Re^m \,\, \mbox{is
onto}, \\
\exists \, X_0 \in \mathcal{S}^n_{\+} \,\, \mbox{such that}\,\,
\mathcal{A}(X_0) = b,\ X_0 \succ \mathbf{0}.
\end{array}
\right.
\end{eqnarray}
From \cite[Theorems $17$ and  $18$]{rockafellar1974conjugate}, we have the following useful lemma.
\begin{lem} \label{lem-innerbonuded}
Suppose that Assumption $\ref{assump-slater}$ holds. Then for any $\alpha \in \Re$, the level set $\mathcal{L}_\alpha :=
  \{ (y,S,Z) \in \Re^m \times \cK^* \times \cP^* \, \mid\, \phi(y,S,Z) \leq \alpha\}$ is a closed and bounded convex set.
\end{lem}

\subsection{A Majorized Semismooth Newton-CG Method }
\def\hC{\widehat{C}}

Consider $(\ty,\tS,\tZ) \in \argmin  \{ \phi(y,S,Z) \,\mid\, y\in\Re^m, S \in \cK^*, Z \in \cP^* \}$.
Let
$$
 \hC = C - \sig^{-1}\tX.
$$
Then we must have $\tZ = \Pi_{\cP^*}(\hC-\cA^*\ty-\tS)$. Therefore, problem (\ref{eq-al-sub-org}) is equivalent to the following optimization problem:
\begin{eqnarray}
   \min \Big\{ \varPhi(y,S) := \inprod{-b}{y} + \frac{\sigma}{2}
     \norm{ \Pi_\cP(\cA^*y + S -\hC)}^2\,\mid\, y\in\Re^m, \; S\in \cK^* \Big\}.
  \label{eq-al-sub}
\end{eqnarray}

In order to introduce our majorized semismooth Newton-CG method for solving (\ref{eq-al-sub}), we need to majorize the second part of  the objective function in (\ref{eq-al-sub}) by a convex, but not necessarily strongly convex,  quadratic function. Specifically, for given $(y^l,S^l) \in \Re^m \times \cK^*$ and $l\ge 0$, since
\begin{eqnarray*}
 \norm{\Pi_{\cP}(\cA^*y+S-\hC)}^2
&\leq& \norm{\Pi_{\cP}(\cA^*y^l+S^l-\hC)}^2 +\;  \norm{\cA^*y+S - \cA^*y^l - S^l}^2
\\[3pt]
& &
  +\;  2 \inprod{ \Pi_{\cP}(\cA^*y^l+S^l-\hC)}{\cA^*y+S- \cA^* y^l -S^l}
  \\[3pt]
  & =&  \norm{\cA^*y+ S + Z^l  -\hC }^2,
\end{eqnarray*}
where  $Z^l: = \Pi_{\cP^*}(\hC-\cA^*y^l-S^l)$, we know that for
$(y,S) \in \Re^m \times \cS^n$,
\begin{eqnarray}
   \varPhi(y,S) &\leq& \inprod{-b}{y}  + \frac{\sigma}{2} \norm{\cA^*y+ S + Z^l -\hC }^2
\nonumber \\[5pt]
  &=& \varPsi_l(y,S) := \inprod{-b}{y} + \frac{\sigma}{2} \norm{\cA^*y+ S + \sig^{-1}\tX-\tC^l}^2,
\label{eq-varPsi}
\end{eqnarray}
where $\tC^l := C - Z^l$.
 Thus $\varPsi_l$ is a majorization function of $\varPhi$ at $(y^l, S^l)$ because $\varPsi_l(y^l,S^l) = \varPhi(y^l,S^l)$ and $\varPsi_l(y,S) \geq \varPhi(y,S)$
$\forall (y,S) \in \Re^m \times \cS^n$.
In order to find an optimal solution for  problem \eqref{eq-al-sub}, for  $l= 0, 1, \ldots$, we solve the following problem
\begin{eqnarray}
   \min \Big\{ \varPsi_l(y,S) \mid y\in\Re^m,\, S \in \cK^* \Big\}.
\end{eqnarray}
Observe that if $(\ty,\tS) \in \argmin  \{ \varPsi_l(y,S) \mid y\in\Re^m,\, S\in \cK^* \}$, then we must have $\tS = \Pi_{\cK^*}(\tC^l-\cA^*\ty-\sig^{-1}\tX)$. Thus we can compute $y^{l+1}$ and $S^{l+1}$ simultaneously as follows:
\begin{subequations} \label{eq-MNCG}
\begin{numcases}{}
   y^{l+1} \in  \argmin \Big\{ \varphi_l(y) :=
      \inprod{-b}{y} + \frac{\sig}{2}\norm{\Pi_{\cK}(\cA^*y+\sig^{-1} \tX -\tC^l)}^2 \mid y \in \Re^m \Big\},  \label{eq-ncg-y}
   \\
   S^{l+1} = \Pi_{\cK^*}(\tC^l-\cA^*y^{l+1}-\sig^{-1}\tX).
\end{numcases}
\end{subequations}
Note that we can only solve  problem (\ref{eq-ncg-y}) inexactly by  an iterative method. Here we will introduce a semismooth Newton-CG (SNCG) method for solving (\ref{eq-ncg-y}). Specifically, for fixed $\tX, \tC \in \cS^n$, we need to consider the following problem of the form
\begin{eqnarray} \label{subprob}
\min \Big\{\varphi(y) :=
   \inprod{-b}{y} + \frac{\sig}{2}\norm{\Pi_{\cK}(\cA^*y+\sig^{-1} \tX-\tC)}^2 \mid y \in \Re^m \Big\}.
\end{eqnarray}
The objective function in (\ref{subprob}) is continuously differentiable and solving (\ref{subprob}) is equivalent to solving the following nonsmooth equation:
\begin{eqnarray} \label{eqn-subnl}
\nabla \varphi(y)= \mathcal{A}\Pi_{\cK} (\tX +
\sigma (\mathcal{A}^* y - \tC)) - b=0, \quad  y \in \Re^m.
\end{eqnarray}
Since $\Pi_{\cK} (\cdot)$ is strongly semismooth \cite{sun2002semismooth}, we can design  a 
SNCG method
as in \cite{zhao2010newton} to
solve  \eqref{eqn-subnl},
 and expect fast superlinear or even quadratic convergence.

 Let $\ty \in \Re^m$ be fixed.
Consider the following eigenvalue decomposition:
\begin{eqnarray}\label{decomp-M}
\tX + \sigma (\mathcal{A}^* \ty - \tC) = Q \, \Gamma_{\ty} \,Q^{\T},
\end{eqnarray}
where $Q\in\R^{n\times n}$ is an orthogonal matrix whose columns are
eigenvectors, and $\Gamma_{\ty}$ is the diagonal matrix of eigenvalues
with the diagonal elements arranged in the nonincreasing order:
$\lambda_1 \geq \lambda_2 \geq
\cdots \geq \lambda_n$. Define the following  index sets
\begin{eqnarray*}
\alpha := \{i \mid \lambda_i>0\}, \quad  \bar{\alpha} :=\{i \mid
\lambda_i\leq 0\}.
\end{eqnarray*} 

\noindent
We define the operator $W^0_{\ty} : \mathcal{S}^n \rightarrow \mathcal{S}^n$ by
\begin{eqnarray} \label{def-W0}
    W^0_{\ty} (H) := Q(\Sigma \circ (Q^{\T}H Q))Q^{\T}, \quad H \in \mathcal{S}^n,
\end{eqnarray}
where $``\circ"$ denotes the Hadamard product of two matrices and
\begin{eqnarray}\label{def-nu} \quad \Sigma = \left[
\begin{array}{cc}
E_{\alpha \alpha}  & \nu_{\alpha \bar{\alpha}}\\[4pt]
\nu^{\T}_{\alpha \bar{\alpha}} & 0
\end{array} \right], \quad
\nu_{ij} := \frac{\lambda_i}{\lambda_i-\lambda_j}, \,\, i \in
\alpha, j \in \bar{\alpha},
\end{eqnarray}
where $E_{\alpha \alpha} \in \mathcal{S}^{|\alpha|}$ is the matrix of ones.
Define
$V^0_{\ty}: \Re^m \rightarrow \mathcal{S}^n$ by
\begin{eqnarray}\label{def-V0}
V^0_{\ty} d := \sigma \mathcal{A} [\, Q (\Sigma \circ (\, Q^{\T}
(\mathcal{A}^*d)\, Q)) \,Q^{\T}], \quad d \in \Re^m.
\end{eqnarray}

 For any $y \in \Re^m$, define
\begin{eqnarray}\label{def-hp2varphi}
\hat{\partial}^2 \varphi(y) := \sigma \mathcal{A} \, \partial
\Pi_{\cK}(\tX +\sigma (\mathcal{A}^* y - \tC))\mathcal{A}^*,
\end{eqnarray}
where $\partial \Pi_{\cK}(\tX +\sigma (\mathcal{A}^* y - \tC))$ is the Clarke subdifferential of $\Pi_{\cK}(\cdot)$ at $\tX +\sigma (\mathcal{A}^* y - \tC)$. Note that from \cite{hiriart1984generalized}, we know that
\begin{eqnarray}
\hat{\partial}^2 \varphi(\ty) h = {\partial}^2 \varphi(\ty) h \quad  \forall \, h \in \Re^m,
\label{eq-Clarke}
\end{eqnarray}
where ${\partial}^2 \varphi(\ty)$ denotes the generalized Hessian of $\varphi$ at $\ty$, i.e., the Clarke subdifferential of $\nabla \varphi$ at $\ty$.
 However, note that \eqref{eq-Clarke} does not mean that $\hat{\partial}^2 \varphi(\ty)  = {\partial}^2 \varphi(\ty)$. Actually, it is unclear to us  whether the latter holds.
 Fortunately, from Pang, Sun, and Sun \cite[Lemma 11]{pang2003semismooth} we know that
$$
W^0_{\ty} \in\partial \Pi_{\cK}(\tX + \sigma (\mathcal{A}^*
\ty - \tC))
$$
and thus
$V^0_{\ty} = \sigma \mathcal{A} W^0_{\ty} \mathcal{A}^* \in
\hat{\partial}^2\varphi (\ty).$

Now we will introduce the  
SNCG algorithm
for solving (\ref{subprob}). Choose $y^0 \in \Re^m$. Then the algorithm can be stated as follows.

\bigskip
\centerline{\fbox{\parbox{\textwidth}{
{\bf Algorithm SNCG}: {\bf A Semismooth Newton-CG Algorithm (SNCG$(y^0, \tX,  \sigma)$).}
\\[5pt]
Given $\mu \in (0, 1/2)$, $\bar{\eta} \in (0, 1)$, $\tau \in (0,1]$,
$\tau_1, \tau_2\in (0,1)$, and $\delta \in (0, 1)$. Perform the $j$th iteration
as follows.
\begin{description}
\item[Step 1.] Given a maximum number of CG iterations $N_j>0$, compute
                    \[ \eta_j := \min(\bar{\eta}, \| \nabla \varphi(y^j)\|^{1+\tau}).\]
                    Apply the conjugate gradient  (CG) algorithm $(CG(\eta_j, N_j))$, to find an approximation solution $d^j$ to
                       \begin{eqnarray}\label{eqn-epsk}
                         (V_j + \varepsilon_j I) \, d= -\nabla \varphi(y^j) ,
                       \end{eqnarray}
                       where $V_j \in \hat{\partial}^2 \varphi(y^j)$ is defined as in (\ref{def-V0}) and $\varepsilon_j := \tau_1 \min\{\tau_2, \|\nabla \varphi
                       (y^j)\|\}$.
\item[Step 2.]  Set $\alpha_j = \delta^{m_j}$, where $m_j$ is the first nonnegative integer $m$ for which
                         \begin{eqnarray}\label{Armijo}
                          \varphi(y^j + \delta^{m} d^j) \leq \varphi(y^j) + \mu \delta^{m}
                           \langle \nabla \varphi(y^j), d^j \rangle.
                          \end{eqnarray}
\item[Step 3.] Set $y^{j+1} = y^j + \alpha_j \, d^j$.
\end{description}
}}}

\bigskip


The convergence results for the above SNCG algorithm are stated
in
Theorems \ref{thm-sub-global} and \ref{thm-sub-rate} below.
We shall omit the proofs as they
can be proved in the same fashion as in \cite[Theorems 3.4 and 3.5]{zhao2010newton}.

\begin{thm}\label{thm-sub-global}
Suppose that Assumption $\ref{assump-slater}$ holds. Then Algorithm SNCG generates a bounded sequence $\{y^j\}$  and any accumulation point $\hat{y}$ of $\{y^j\}$ is an optimal solution to  problem {\rm (\ref{subprob})}.
\end{thm}
\unskip

\begin{thm}\label{thm-sub-rate}
Suppose that Assumption $\ref{assump-slater}$ holds. Let $\hat{y}$ be an accumulation point of the infinite sequence $\{y^j\}$ generated by Algorithm SNCG for solving the problem {\rm (\ref{subprob})}. Suppose that at each step $j\ge 0$, when the   CG  algorithm terminates, the tolerance $\eta_j$ is achieved (e.g., when $N_j = m+1$), i.e.,
\begin{eqnarray}\label{tol-eta}
\|\nabla \varphi(y^{j}) + (V_j + \varepsilon_j I) \, d^j \| \leq
\eta_j.
\end{eqnarray}
Assume that the constraint nondegenerate condition
\begin{eqnarray}
\mathcal{A} \, {\rm lin}(\mathcal{T}_{\cK}
(\widehat{W})) = \Re^m
\end{eqnarray}
holds at $\widehat{W}:= \Pi_{\cK} (\tX + \sigma (\mathcal{A}^* \hy - \tC))$, where ${\rm lin} (\mathcal{T}_{\cK} (\widehat{W}))$ denotes the lineality space of  $\mathcal{T}_{\cK} (\widehat{W})$. Then the whole sequence $\{y^j\}$ converges to $\hat{y}$ and
\begin{eqnarray}
\|y^{j+1} - \hat{y} \| = O( \|y^{j} - \hat{y} \|^{1+\tau}).
\end{eqnarray}
\end{thm}

\bigskip
Given $\xi_1 \in (0,1)$ and $\xi_2 \in (0,\infty)$, we will use the following stopping criteria  for terminating Algorithm SNCG:
\begin{itemize}
\item[(A1)]  $\varphi_l(y^{l+1}) \leq \varphi_l(y^l) - \frac{\xi_1}{2} \abs{\inprod{\nabla \varphi_l(y^l)}{y^{l+1}-y^l}}$,

\item[(A2)] $\norm{\nabla \varphi_{l}(y^{l+1})} \leq \xi_2 (\varphi_l(y^l) - \varphi_l(y^{l+1}) )^\frac{1}{2}$.

\end{itemize}

We can now state our majorized semismooth Newton-CG method for solving (\ref{eq-al-sub})  as follows:

\bigskip
\centerline{\fbox{\parbox{\textwidth}{
{\bf Algorithm MSNCG}: {\bf A Majorized Semismooth Newton-CG Algorithm (MSNCG$(y^0,S^0,Z^0,\tX,\sig)$).}
\\[5pt]
 Given  $\xi_1 \in (0,1)$, $\xi_2 \in (0,+\infty)$. Perform the $l$th iteration as
follows.
\begin{description}
\item[Step 1.] Starting with $y^l$ as the initial point, apply Algorithm SNCG to minimize $\varphi_l(\cdot)$ to find
   $y^{l+1}={\rm SNCG}(y^l,\tX,\sigma)$  satisfying (A1) and (A2).
\item[Step 2.] Compute
 $S^{l+1}  := \Pi_{\cK^*}(C-\cA^*y^{l+1}-Z^l-\sig^{-1}\tX)$ and $Z^{l+1} := \Pi_{\cP^*}(C-\cA^*y^{l+1}-S^{l+1}-\sig^{-1} \tX)$.
\end{description}
}}}
\bigskip

Next, we establish the convergence of Algorithm MSNCG.
For notational convenience, for any $y \in \Re^m$ and $S \in \cS^n$, let $Y: = (y,S)$.  Define the linear map $\cM : \Re^m \times \cS^n \to \cS^n$ by
\begin{eqnarray}
\cM Y := \cA^*y + S \quad \forall \, Y = (y,S) \in \Re^m \times \cS^n.
\end{eqnarray}
Let $B := (b,0) \in \Re^m \times \cS^n$ and $\cC := \Re^m \times \cK^*$.  Then problem (\ref{eq-al-sub}) is equivalent to
\begin{eqnarray}
\min \Big\{ \varPhi(Y) := \inprod{-B}{Y} + \frac{\sig}{2}\norm{\Pi_\cP(\cM Y + \sig^{-1}\tX -C)}^2 \,\mid\,  Y \in \cC \Big\}
  \label{eq-al-sub-mat}
\end{eqnarray}
and  the function $\varPsi_l(y,S)$ in \eqref{eq-varPsi} can be rewritten as
\begin{eqnarray*}
\varPsi_l(Y) =\varPsi_l(y,S) = -\inprod{B}{Y} + \frac{\sig}{2}\norm{\cM Y +\sig^{-1}\tX - C + Z^l}^2.
\end{eqnarray*}
Furthermore, $\hY = (\hy,\hS)$ is an optimal solution of
\begin{eqnarray}
\min \Big\{ \varPsi_l(Y)  \,\mid\, Y \in \cC  \Big\}
  \label{eq-al-sub-mat-major}
\end{eqnarray}
if and only if $\hy$ is an optimal solution of problem (\ref{eq-ncg-y}) and $\hS = \Pi_{\cK^*}(\tC^l-\cA^*\hy-\sig^{-1}\tX).$
In addition, conditions (A1) and (A2) are equivalent to the following two conditions, respectively,
\begin{eqnarray}
  && \varPsi_l(Y^{l+1}) \leq \varPsi_l(Y^{l}) - \frac{\xi_1}{2} \abs{\inprod{\nabla \varPsi_l(Y^{l})}{Y^{l+1}-Y^l}},
  \label{eq-Theta-dec}
  \\[3pt]
  && \norm{Y^{l+1} - \Pi_\cC(Y^{l+1} - \nabla \varPsi_l(Y^{l+1}))} \leq \xi_2 (\varPsi_l(Y^{l}) - \varPsi_l(Y^{l+1}) )^\frac{1}{2}. \label{eq-GradTheta-dec}
\end{eqnarray}

\begin{lem} \label{lem-MNCG-criteria}
Suppose that Assumption $\ref{assump-slater}$ holds. Then for {Algorithm MSNCG}, {\rm (A1)} and {\rm (A2)} are achievable.
\end{lem}
\begin{proof}
 If $Y^{l} - \Pi_\cC(Y^{l} - \nabla \varPsi_l(Y^{l})) = 0$, then one can take $Y^{l+1} = Y^l$ to satisfy (A1) and (A2).

Next, we assume that $Y^{l} - \Pi_\cC(Y^{l} - \nabla \varPsi_l(Y^{l})) \neq 0$. Then $Y^l$ is not an optimal solution of problem (\ref{eq-al-sub-mat-major}).
Let $\hY$ be an arbitrary optimal solution of problem (\ref{eq-al-sub-mat-major}).  Then $\hY = \Pi_\cC(\hY - \nabla \varPsi_l(\hY))$. So
$
\inprod{Y^l - \hY}{(\hY - \nabla \varPsi_l(\hY))-\hY} \leq 0,
$
i.e.,
\begin{eqnarray*}
 \inprod{\nabla \varPsi_l(\hY)}{Y^l - \hY} \geq 0,
\end{eqnarray*}
which implies
\begin{eqnarray}
 \inprod{\nabla \varPsi_l(Y^l)}{Y^l - \hY} \geq  \inprod{\nabla \varPsi_l(Y^l) -\nabla\varPsi_l(\hY)}{Y^l - \hY}
=
\sig \norm{\cM(\hY-Y^l)}^2.
\label{eq-gradTheta}
\end{eqnarray}
Since
\begin{eqnarray} \label{eq-ThetaYhat}
 \varPsi_l(Y^l) > \varPsi_l(\hY) = \varPsi_l(Y^l) + \inprod{\nabla \varPsi_l(Y^l)}{\hY - Y^l} + \frac{\sig}{2}\norm{\cM(\hY-Y^l)}^2,
\end{eqnarray}
we obtain that
$
\inprod{\nabla \varPsi_l(Y^l)}{\hY - Y^l} + \frac{\sig}{2}\norm{\cM(\hY-Y^l)}^2 < 0.
$
This implies
\begin{eqnarray} \label{eq-gradTheta-neg}
\inprod{\nabla \varPsi_l(Y^l)}{\hY - Y^l} < 0.
\end{eqnarray}
Then by using (\ref{eq-gradTheta}), (\ref{eq-ThetaYhat}), (\ref{eq-gradTheta-neg}) and the fact that $\hY = \Pi_\cC(\hY - \nabla \varPsi_l(\hY))$, we know that for given $\xi_1 \in (0,1)$ and $\xi_2 \in (0,\infty)$, there exists $\delta > 0$ such that
\begin{eqnarray*}
&& \varPsi_l(Y) \leq \varPsi_l(Y^l) - \frac{\xi_1}{2} \abs{\inprod{\nabla \varPsi_l(Y^l)}{Y - Y^l}},
\\
&& \norm{Y - \Pi_\cC(Y - \nabla \varPsi_l(Y))} \leq \xi_2 (\varPsi_l(Y^{l}) - \varPsi_l(Y) )^\frac{1}{2},
\end{eqnarray*}
for all $Y \in \cC$ satisfying  $\norm{Y - \hY} < \delta$.
Let $\{\ty^j\}_{j=0}^{+\infty}$ be the sequence generated by SNCG$(\ty^0,\tX,\sigma)$ with $\ty^0 := y^l$. For each $j\ge 0$, let $\tS^j := \Pi_{\cK^*}(C-\cA^*\ty^{j}-Z^l-\sig^{-1}\tX)$. Then by Theorem \ref{thm-sub-global}, we know that $\{(\ty^j,\tS^j)\}$ is a bounded sequence and any accumulation point of $\{(\ty^j,\tS^j)\}$,  say $\hY := (\hy,\hS)$,   is an optimal solution to problem {\rm (\ref{eq-al-sub-mat-major})}. Thus there exists a sufficiently large $j$ such that $Y^{l+1} := (\ty^j,\tS^j)$  satisfying (\ref{eq-Theta-dec}) and (\ref{eq-GradTheta-dec}).
\end{proof}

\begin{thm} \label{thm-MNCG}
Suppose that Assumption $\ref{assump-slater}$ holds. Let Algorithm MSNCG be executed with stopping criteria {\rm (A1)} and {\rm (A2)}.  Then it generates a bounded sequence $\{(y^l,S^l,Z^l)\}$ and any accumulation point $(\hy,\hS)$ of $\{(y^l,S^l)\}$ is an optimal solution to problem {\rm (\ref{eq-al-sub})} and hence $(\hy,\hS,\hZ)$ is an optimal solution to problem {\rm (\ref{eq-al-sub-org})}, where $\hZ := \Pi_{\cP^*}(C-\cA^*\hy-\hS-\sig^{-1} \tX)$. Furthermore, $\norm{Z^{l+1}-Z^l} \to 0$ as $l \to +\infty$.
\end{thm}

\begin{proof}
By (\ref{eq-Theta-dec}), we have $\varPhi(Y^{l+1}) \leq \varPsi_l(Y^{l+1}) \leq \varPsi_l(Y^{l}) - \frac{\xi_1}{2} \abs{\inprod{\nabla \varPsi_l(Y^{l})}{Y^{l+1}-Y^l}} = \varPhi(Y^{l}) - \frac{\xi_1}{2} \abs{\inprod{\nabla \varPsi_l(Y^{l})}{Y^{l+1}-Y^l}}$. Hence, the sequence $\{\varPhi(Y^{l})\}$ is nonincreasing.

By Lemma \ref{lem-innerbonuded}, we know that the level set $\mathcal{L} :=  \{ Y \in \cC \, \mid\, \varPhi(Y) \leq \varPhi(Y^0)\}$ is a   closed and bounded convex set. Then the sequence $\{Y^l\}$ is bounded and so is the sequence $\{Z^l\}$. Let $\hY$ be an accumulation point of $\{Y^l\}$.
Then $\varPhi(Y^l) \to \varPhi(\hY)$ and $\inprod{\nabla \varPsi_l(Y^{l})}{Y^{l+1}-Y^l} \to 0$ as $l \to \infty$. Furthermore, $\varPsi_l(Y^{l}) - \varPsi_l(Y^{l+1}) \rightarrow 0$ as $l \to \infty$.

By noting that
\begin{eqnarray} \label{eq-TY-lp1}
  \varPsi_l(Y^{l+1}) = \varPsi_l(Y^l) + \inprod{\nabla \varPsi_l(Y^{l})}{Y^{l+1}-Y^l} + \frac{\sig}{2} \norm{\cM (Y^{l+1}-Y^l)}^2,
\end{eqnarray}
we get from (\ref{eq-Theta-dec}) that
\begin{eqnarray} \label{eq-MY-neg}
     \inprod{\nabla \varPsi_l(Y^{l})}{Y^{l+1}-Y^l} + \frac{\sig}{2} \norm{\cM (Y^{l+1}-Y^l)}^2 \leq - \frac{\xi_1}{2} \abs{\inprod{\nabla \varPsi_l(Y^{l})}{Y^{l+1}-Y^l}} \leq 0.
\end{eqnarray}
Since $\inprod{\nabla \varPsi_l(Y^{l})}{Y^{l+1}-Y^l} \to 0$ as $l \to \infty$, we  obtain from (\ref{eq-MY-neg}) that
  \begin{eqnarray} \label{eq-MYtozero}
   \norm{\cM (Y^{l+1}-Y^l)}  \to 0 \quad  {\rm as} \quad l \to \infty.
  \end{eqnarray}
For any  $l\ge 0$, denote $\Delta_{l}: = Y^{l} - \Pi_\cC(Y^{l} - \nabla \varPhi(Y^{l}))$. Then we have
  \begin{eqnarray*}
   \norm{\Delta_{l+1}} & \leq & \norm{Y^{l+1} - \Pi_\cC(Y^{l+1} - \nabla \varPsi_l(Y^{l+1}))}
   \\
   &  & +\; \norm{\Pi_\cC(Y^{l+1} - \nabla \varPsi_l(Y^{l+1}))-\Pi_\cC(Y^{l+1} - \nabla \varPhi(Y^{l+1}))}
   \\
   & \leq & \xi_2 (\varPsi_l(Y^{l}) - \varPsi_l(Y^{l+1}) )^\frac{1}{2} + \norm{\nabla \varPsi_l(Y^{l+1}) - \nabla \varPhi(Y^{l+1})}.
  \end{eqnarray*}
By direct computations, we have for $l\ge 1$,
  \begin{eqnarray*}
    && \hspace{-0.7cm} \norm{\nabla \varPsi_l(Y^{l+1}) - \nabla \varPhi(Y^{l+1})}
    \\[3pt]
    &=&  \norm{\sig \cM^*(\cM Y^{l+1}+\sig^{-1}\tX-C+Z^l) -
\sig \cM^*(\Pi_{\cP}(\cM Y^{l+1}+\sig^{-1} \tX-C))}
    \\[3pt]
& = &
 \norm{\sig\cM^*( Z^l- Z^{l+1})}
    \; \leq \; \sig \norm{\cM^*} \norm{\cM (Y^{l+1}-Y^l)},
  \end{eqnarray*}
  where we have used the fact that
\begin{eqnarray}\label{eq:Zdifference}
\norm{Z^{l+1}-Z^l} &=& \norm{\Pi_{\cP^*}(C- \sigma^{-1}\tX-\cM Y^{l+1} ) - \Pi_{\cP^*}(C- \sigma^{-1}\tX-\cM Y^{l} )} \nonumber
\\
&\leq& \norm{\cM(Y^{l+1}-Y^l)}.
\end{eqnarray}
Thus, by (\ref{eq-MYtozero}) and the fact that $(\varPsi_l(Y^{l}) - \varPsi_l(Y^{l+1}) )^\frac{1}{2} \rightarrow 0$ as $l \to \infty$, we derive that $\norm{\Delta_{l+1}} \to 0$ as $l \to \infty$. Since $\hY$ is an accumulation point of $\{Y^l\}$, we obtain that $\hY - \Pi_\cC(\hY - \nabla \varPhi(\hY)) = 0$. By the convexity of $\varPhi$, $\hY$ is an optimal solution of problem (\ref{eq-al-sub-mat}).

Finally, by using (\ref{eq-MYtozero}) and (\ref{eq:Zdifference}), we know that $\norm{Z^{l+1}-Z^l} \to 0$ as $l \to \infty$.
\end{proof}

\section{A Majorized Semismooth Newton-CG Augmented Lagrangian Method }
\label{sec-SDPNALP}

For any $k \geq 0$ and $ (y, S,Z) \in\Re^m\times  \cS^n\times  \cS^n$, denote
\begin{eqnarray}
   \phi_k(y,S,Z) &:=& L_{\sig_k}(y,S,Z;X^k),
   \\
   \hat{\phi}_k(y,S,Z) &:=& \begin{cases}
    L_{\sig_k}(y,S,Z; X^k) & {\rm if } \; (y,S,Z) \in \varOmega := \Re^m \times \cK^* \times \cP^*,
    \\
    +\infty & {\rm otherwise}.
   \end{cases}
\end{eqnarray}
Since the inner problems in (\ref{eq-al-sub-org}) are solved inexactly, we will use the following standard
stopping criteria considered in \cite{rockafellar1976monotone,rockafellar1976augmented} to terminate Algorithm MSNCG:

\begin{itemize}

\item[(B1)]  $\hat{\phi}_k(y^{k+1},S^{k+1},Z^{k+1}) - \inf \hat{\phi}_k  \leq \epsilon_k^2$, $\epsilon_k \geq 0$, \vspace*{2pt} $\sum^{\infty}_{k=0} \epsilon_k < \infty$.

\item[(B2)]  $\hat{\phi}_k(y^{k+1},S^{k+1},Z^{k+1}) - \inf \hat{\phi}_k  \leq (\delta_k^2/2\sigma_k)  \|X^{k+1} -
X^k\|$, $\delta_k \geq 0$, \vspace*{2pt} $\sum^{\infty}_{k=0} \delta_k < \infty$.

\item[(B3)] ${\rm dist}(0,\partial \hat{\phi}_k(y^{k+1},S^{k+1},Z^{k+1}))  \leq (\delta_k^{'}/\sigma_k)  \|X^{k+1} -
X^k\|$, $0 \leq \delta^{'}_k \rightarrow 0$.
\end{itemize}

Just like SDPNAL, each iteration of the MSNCG algorithm can be quite expensive.
Thus it is crucial for us
to find a reasonably good initial point to warm start Algorithm SDPNAL$+$. We can certainly do so by  solving the inner problem (\ref{eq-al-sub-org}) by using any  gradient descent type method. However,  for this purpose we find that ADMM+ introduced by Sun, Toh and Yang \cite{sun2014convergent} is usually more efficient than other choices.
Now we can present our SDPNAL$+$  algorithm as follows.

\bigskip
\centerline{\fbox{\parbox{\textwidth}{
{\bf Algorithm SDPNAL$+$}: {\bf A Majorized Semismooth Newton-CG Augmented Lagrangian Algorithm (SDPNAL$+$$(y^0,S^0,Z^0,X^0,\sig_0)$)}
\begin{description}
\item[Stage {\rm 1.}]
Use ADMM+ to generate an initial point
\\[5pt]
\centerline{
$(y^{0},S^{0},Z^0,X^{0},\sig_0) \leftarrow  {\rm ADMM+} (y^{0},S^{0},Z^0,X^{0},\sig_0).$}
\item[Stage {\rm 2.}] For $k=0,\ldots,$ perform the $k$th iteration as follows:
\begin{enumerate}
\item[(a)] Using $(y^k,S^k,Z^k)$ as the initial point, apply Algorithm MSNCG to minimize $\hat{\phi}_k(\cdot)$ to find
   $(y^{k+1},S^{k+1},Z^{k+1})=$ MSNCG $(y^k,S^k,Z^k,X^k,\sigma_k)$ and $X^{k+1} = X^k + \sig_k(\cA^*y^{k+1}+ S^{k+1}+Z^{k+1}-C)$ satisfying (B1), (B2) or (B3).
\item[(b)] Update $\sigma_{k+1} =\rho \sigma_k$ for some $\rho>1$ or $\sigma_{k+1}=\sigma_k$. 
\end{enumerate}
\end{description}
}}}

\bigskip
	
\begin{remark}
As mentioned in the introduction, if {\rm (P)} is  reformulated as a standard SDP and Algorithm SDPNAL$+$ is applied to this reformulated form, then SDPNAL$+$ reduces to SDPNAL proposed in \cite{zhao2010newton}.
\end{remark}

We can obtain   similar theorems on the convergence of SDPNAL$+$ as  SDPNAL (\cite[Theorems 4.1 and 4.2]{zhao2010newton}). The global convergence of Algorithm SDPNAL$+$ follows from Rockafellar \cite[Theorem $1$]{rockafellar1976monotone} and  \cite[Theorem $4$]{rockafellar1976augmented} without much difficulty.
\begin{thm} \label{thm-SDPNALP-convergence}
Suppose that Assumption $\ref{assump-slater}$ holds. Let Algorithm SDPNAL$+$ be executed with stopping criterion {\rm (B1)}.
If there exists $(y_0,S_0,Z_0) \in \Re^m \times \cS^n_{+} \times \cS^n$ such that
\begin{eqnarray} \label{slater-dual}
\cA^*(y_0) + S_0 + Z_0 = C, \; S_0 \succ 0, \; Z_0 \in relint(\cP^*),
\end{eqnarray}
then the sequence $\{X^k\} \subset \cP$ generated
by Algorithm SDPNAL$+$ is bounded and $\{ X^k \}$ converges to
$\overline{X}$, where $\overline{X}$ is some optimal solution to
{\rm (P)}, and $\{(y^k,S^k,Z^k)\}$ is asymptotically minimizing for
{\rm (D)} with $\max({\rm P})=\inf({\rm D})$.

If $\{X^k\}$ is bounded, then the sequence $\{(y^k,S^k,Z^k)\}$ is also bounded, and all of its accumulation points of the sequence $\{(y^k,S^k,Z^k)\}$ are optimal solutions to {\rm (D)}.
\end{thm}

Next we state the local linear convergence of Algorithm SDPNAL$+$.

\begin{thm} \label{thm-SDPNALP-global-rate}
Suppose that Assumption $\ref{assump-slater}$ holds. Let Algorithm SDPNAL$+$ be executed with stopping criteria {\rm (B1)} and
{\rm (B2)}. Assume that {\rm (D)} satisfies condition {\rm (\ref{slater-dual})}. If the second order sufficient conditions (in the sense of the conditions in {\rm \cite[Theorem $3.137$]{bonnans2000perturbation}}) holds at $\overline{X}$, where
$\overline{X}$ is an optimal solution to {\rm (P)}, then the
generated sequence $\{X^k\} \subset\cP$ is bounded
and $\{ X^k \}$ converges to the unique optimal solution $\overline{X}$
with $\max({\rm P}) = \min({\rm D})$, and
$$
\|X^{k+1} - \overline{X}\| \leq \theta_{\infty} \|X^k - \overline{X} \|
\quad \mbox{$\forall\ k$ sufficiently large,}
$$
for some $\theta_{\infty} \in [0,1)$ with the property that $\theta_\infty \ll 1$ if  $\sigma_k\to \sigma_\infty $ for any sufficiently large $\sigma_\infty$. The conclusions of Theorem {\rm \ref{thm-SDPNALP-convergence}} about $\{(y^k,S^k,Z^k)\}$ are also valid.
\end{thm}
\unskip

\begin{proof}
The conclusions of Theorem \ref{thm-SDPNALP-global-rate} follow
from the results in \cite[Theorem 2]{rockafellar1976monotone} and
\cite[Theorem 5 and Proposition 3]{rockafellar1976augmented} combined with \cite[Theorem $3.137$]{bonnans2000perturbation}.
\end{proof}

\section{Numerical Experiments}
\label{sec-numerical}

\subsection{SDP$+$ and SDP Problem Sets}

In our numerical experiments, we test the following SDP$+$ and SDP problem sets.

(i) SDP$+$ problems coming from
the relaxation of a binary integer nonconvex quadratic (BIQ) programming:
\begin{eqnarray}
\min\Big\{ \frac{1}{2}x^TQx + \inprod{c}{x} \mid  x\in\{0,1\}^{n-1}\Big\}.
\end{eqnarray}
This problem has  been shown in \cite{burer2009copositive}  that under some mild
assumptions, it can equivalently be reformulated as the following
completely positive programming (CPP) problem:
\begin{eqnarray}
\min\Big\{ \frac{1}{2} \inprod{Q}{X_0} + \inprod{c}{x} \mid
 \diag(X_0) = x, X = [X_0, x; x^T, 1]\in
{\cal C}_{pp}^n\Big\},
\end{eqnarray}
where
${\cal C}_{pp}^n$ denotes the $n$-dimensional completely positive cone.
It is well known that even though ${\cal C}_{pp}^n$ is convex, it is
computationally intractable. To solve the CPP problem,  one would typically
relax ${\cal C}_{pp}^n$ to $\S^n_+\cap \cN$, and the relaxed problem
has the form (P):
\begin{eqnarray}
\begin{array}{ll}
   \min & \frac{1}{2} \inprod{Q}{X_0} + \inprod{c}{x} \\[5pt]
  {\rm s.t.} &  {\rm diag}(X_0) - x = 0, \; \alpha = 1,
  \quad X = \left[ \begin{array}{cc} X_0 & x \\ x^T & \alpha  \end{array}\right] \in \cS^n_+,\quad
  X \in \cP,
  \end{array}
  \label{eq-BIQ-SDP}
\end{eqnarray}
where the polyhedral cone $\cP = \{X \in \cS^n \mid X \geq 0 \}$. In our numerical experiments, the test data for $Q$ and $c$ are taken from Biq Mac Library  maintained by Wiegele, which is available at {\tt http://biqmac.uni-klu.ac.at/biqmaclib.html}.

\bigskip
(ii) SDP and SDP$+$ problems arising from the relaxation of maximum stable set problems. Given a graph $G$ with edge set $\mathcal{E}$, the SDP and SDP$+$ relaxation $\theta(G)$ and $\theta_+(G)$ of the maximum stable set problem are given by
\begin{eqnarray}
\theta(G)  &=& \max \{ \inprod{ee^T}{X} \,:\ \inprod{E_{ij}}{X} =0,
  (i,j) \in \mathcal{E}, \,\inprod{I}{X} =1, \, X \in \cS^n_+ \},
   \label{eq-theta}
   \\[5pt]
 \theta_+(G)  &=& \max \{ \inprod{ee^T}{X} \,:\ \inprod{E_{ij}}{X} =0,
  (i,j) \in \mathcal{E}, \,\inprod{I}{X} =1, \, X \in \cS^n_+,\, X \in \cP \},\quad
  \label{eq-theta-plus}
\end{eqnarray}
where $E_{ij}=e_ie_j^{\T} + e_je_i^{\T}$ and $e_i$ denotes the $i$th column
 of the identity matrix, $\cP = \{X \in \cS^n \mid X \geq 0 \}$. In our numerical experiments, we test the graph instances $G$ considered in
 \cite{Sloane}, \cite{toh2004}, and \cite{DIMACS}.

 \bigskip
(iii) SDP$+$ relaxation for computing lower
bounds for quadratic assignment problems (QAPs). Let $\Pi$ be the set of $n\times n$ permutation matrices. Given
matrices $A, B \in \S^n$, the QAP is given by
\begin{eqnarray}
 v^*_{\rm QAP} := \min \{ \langle X, AXB \rangle \, :\, X \in \Pi\}.
 \label{eq-QAP}
\end{eqnarray}
For a matrix $X = [x_1,\dots,x_n]\in \Re^{n\times n}$, we will
identify it with the $n^2$-vector $x = [x_1; \dots; x_n]$. For a
matrix $Y\in R^{n^2 \times n^2}$, we let $Y^{ij}$ be the $n\times n$
block corresponding to $x_i x_j^T$ in the matrix $xx^T$. It is shown
in \cite{povh2009copositive} that $v^*_{\rm QAP}$ is bounded below by the
following number generated from the SDP$+$ relaxation of (\ref{eq-QAP}):
\begin{eqnarray}
  \begin{array}{rl}
  v \,:=\, \min & \inprod{B \otimes A}{Y} \\[5pt]
  {\rm s.t.} &  \sum_{i=1}^n Y^{ii} = I, \,
  \inprod{I}{Y^{ij}} = \delta_{ij} \quad \forall \, 1\leq i \leq j\leq n, \\[5pt]
  &
  \inprod{E}{Y^{ij}} = 1\quad \forall\, 1\leq i \leq j\leq n,
  \quad Y \in \cS^n_+,\; Y \in \cP,
  \end{array}
  \label{eq-PR}
\end{eqnarray}
where $E$ is the matrix of ones, and $\delta_{ij} = 1$ if $i=j$,
and $0$ otherwise, $\cP = \{X \in \cS^{n^2} \mid X \geq 0 \}$. In our numerical experiments, the test instances $(A,B)$ are taken
from the QAP Library \cite{QAPLIB}.

\bigskip
(iv) SDP$+$ relaxations of  clustering problems (RCPs)  described in \cite[eq. (13), up to a constant]{peng2007approximating}:
\begin{eqnarray}
\min \Big\{ \inprod{-W}{X} \mid Xe = e, \inprod{I}{X} = K, X \in \S^n_+,
X \in \cP \Big\},
\end{eqnarray}
where $W$ is the so-called affinity matrix whose entries represent the similarities of the objects in the dataset, $e$ is the vector of ones, and $K$ is the number of clusters, $\cP = \{X \in \cS^n \mid X \geq 0 \}$.  All the data sets we tested are from the UCI Machine Learning Repository (available at {\tt http://archive.ics.uci.edu/ml/datasets.html}). For some large data instances, we only select the first $n$ rows. For example, the original data instance ``{\tt spambase}" has 4601 rows, we select the first 1500 rows to obtain the test problem ``{\tt spambase-large.2}"
for which the number ``{\tt 2}" means that there are $K = 2$ clusters.

\bigskip
(v) SDP+ problems arising from the SDP relaxation of frequency assignment problems (FAPs) \cite{eisenblatter2002frequency}. The explicit description of the SDP in the form (\DD) is given in \cite[eq. (5)]{burer2003computational}:
\begin{eqnarray} \label{eq-fap}
\begin{array}{ll}
   \max & \inprod{(\frac{k-1}{2k})L(G,W)-\frac{1}{2}{\rm Diag}(We)}{ X} \\[5pt]
  {\rm s.t.} & {\rm diag}(X) = e, X \in \cS^n_+,
  \\[5pt]
  & -E^{ij} \bullet X = 2/(k-1) \quad \forall (i,j) \in U \subseteq E,
  \\[5pt]
  & -E^{ij} \bullet X \leq 2/(k-1) \quad \forall (i,j) \in E \setminus U,
  \end{array}
\end{eqnarray}
where $k>1$ is an integer, $L(G,W) := {\rm Diag}(We)-W$ is the Laplacian matrix, $E^{ij}=e_ie_j^T + e_je_i^T$ with $e_i \in \Re^n$ being the $i$th standard unit vector and $e \in \Re^n$ is the vector of all ones. Let
\begin{eqnarray}
M_{ij} =
\begin{cases}
-\frac{1}{k-1} &  \forall (i,j) \in E, \\
0 & {\rm otherwise}.
\end{cases}
\end{eqnarray}
Then (\ref{eq-fap}) is equivalent to
\begin{eqnarray} \label{eq-fap-sdpp}
 \begin{array}{ll}
    \max & \inprod{(\frac{k-1}{2k})L(G,W)-\frac{1}{2}{\rm Diag}(We)}{ X}
           \\[5pt]
      {\rm s.t.} & {\rm diag}(X) = e, X  \in \cS^n_+, X - M \in \cP,
 \end{array}
\end{eqnarray}
where $\cP = \{X \in \cS^n \mid X_{ij} = 0, \forall (i,j) \in U; X_{ij} \geq  0, \forall (i,j) \in E \setminus U \} $.

We should mention that we can easily extend our algorithm to handle the following  more general SDP$+$ problem:
\begin{eqnarray}
   \min \Big\{ \inprod{C}{X} \mid \cA(X) = b,\; X  \in \cK, \; X - M \in \cP \Big\},
  \label{eq-primal-general}
\end{eqnarray}
where  $M \in \cX$ is a given matrix. Thus \eqref{eq-fap-sdpp} can also be solved by our proposed  algorithm.

\bigskip
(vi) SDP relaxations for rank-1 tensor approximations (R1TA) \cite{nie2013semidefinite}:
\begin{eqnarray} \label{eq-R1TA}
\max \Big\{ \inprod{f}{y} \mid  M(y) \in \cS^n_+ , \; \inprod{g}{y} = 1 \Big\},
\end{eqnarray}
where $y \in \Re^{\mathbb{N}_m^n}$, $M(y)$ is a linear pencil in $y$.  The dual of (\ref{eq-R1TA}) is given by
\begin{eqnarray} \label{eq-R1TA-dual}
\min \Big\{ \gamma \mid \gamma g - f = M^*(X), \; X \in \cS^n_+ \Big\}.
\end{eqnarray}
It is shown in \cite{nie2012regularization} that  (\ref{eq-R1TA-dual}) can be transformed into a standard SDP (up to a constant):
\begin{eqnarray}
  \quad  \min \Big\{ \inprod{C}{X}  \mid \cA(X) = b,\; X \in \cS^n_+ \Big\},
\end{eqnarray}
where $C$ is a constant matrix and $\cA$ is a linear map, which  depend on $M, f, g$.

\subsection{Numerical Results}

In this subsection, we compare the performance of our SDPNAL$+$ algorithm with
two other  competitive
publicly available  first order methods based codes  for solving large-scale SDP$+$ and SDP problems:  an ADMM based solver, called SDPAD\footnote{\url{http://www.bicmr.org/~wenzw/code/SDPAD-release-beta2.zip}} (release-beta2, released in December 2012) developed in \cite{wen2010alternating} and a two-easy-block-decomposition hybrid proximal extragradient method, which was called 2EBD-HPE\footnote{\url{www2.isye.gatech.edu/~cod3/CamiloOrtiz/Software_files/2EBD-HPE_v0.2/2EBD-HPE_v0.2.zip}} (v0.2, released on May 31, 2013) and we call it 2EBD  here, introduced in \cite{monteiro2013first}.
Since  we use the  convergent ADMM with $3$-block constraints  introduced by Sun et al. \cite{sun2014convergent} (which was called ADMM3c but we call it ADMM$+$ here to indicate that it is an enhanced version of ADMM with convergence guantantee) to warm start SDPNAL$+$, we also list the numerical results obtained by running ADMM$+$ alone for the purpose of demonstrating the power and the importance of the proposed majorized semismooth Newton-CG algorithm for solving difficult SDP$+$ and SDP problems.

 All our computational results for   the tested SDP$+$ and SDP problems are obtained by running {\sc Matlab} on a Linux server  having 6 cores with 12 Intel Xeon  X5650 processors at 2.67GHz and 32G RAM.

Note that numerically it is difficult to compute ${\rm dist}(0,\partial \hat{\phi}_k(W^{k+1}))$ in the criterion (B3) for terminating Algorithm MSNCG  directly,
where $W^{k+1}=(y^{k+1},S^{k+1},Z^{k+1})$. Fortunately, by using the fact that for any closed convex cone $\cC \subseteq \cS^n$ and $X \in \cC$, it holds that  $\cC \subseteq \cT_{\cC}(X)$, we have from \cite{rockafellar1976augmented} that
\begin{eqnarray}
& & \hspace{-0.7cm} {\rm dist}(0,\partial \hat{\phi}_k(W^{k+1}))
\;=\; \norm{\Pi_{\cT_\varOmega(W^{k+1})}(-\nabla \phi_k(W^{k+1}))}  \nonumber
\\[5pt]
&\leq& \norm{\Pi_{\cT_{\Re^m}(y^{k+1})}(-\nabla_y \phi_k(W^{k+1}))} + \norm{\Pi_{\cT_{\cK^*}(S^{k+1})}(-\nabla_{S} \phi_k(W^{k+1}))}  + \norm{\Pi_{\cT_{\cP^*}(Z^{k+1})}(-\nabla_{Z} \phi_k(W^{k+1}))}  \nonumber
\\[5pt]
&\leq& \norm{-\nabla_y \phi_k(W^{k+1})} + \norm{\Pi_{\cK^*}(-\nabla_{S} \phi_k(W^{k+1}))}    + \norm{\Pi_{\cP^*}(-\nabla_{Z} \phi_k(W^{k+1}))}  \nonumber
\\[5pt]
&=& \norm{ \cA ( X^k+ \sig_k R_D^{k+1})-b}
 + \sig_k \norm{\Pi_{\cK^*}(-R_D^{k+1} - \sig_k^{-1}X^k)}. \label{eq-dist-grad}
\end{eqnarray}
where $R_D^{k+1} = \cA^* y^{k+1} + S^{k+1} + Z^{k+1}-C$.
In order to avoid computing $\Pi_{\cK^*}(-R_D^{k+1} - \sig_k^{-1}X^k)$, we majorize the second part of (\ref{eq-dist-grad})  by a  simpler term. Specifically, by using the fact that for any  $X \in \cK$ and $Y \in \cS^n$, $\norm{\Pi_{\cK^*}(Y-X)} = \norm{\Pi_{\cK^*}(Y-X) - \Pi_{\cK^*}(-X)} \leq \norm{Y}$, we have
\begin{eqnarray*}
&&\hspace{-0.7cm}
  \norm{\Pi_{\cK^*}(-R_D^{k+1} - \sig_k^{-1}X^k)}
=  \norm{\Pi_{\cK^*}((\tZ^k - Z^{k+1}) - ( \cA^* y^{k+1} + S^{k+1} +\tZ^{k} + \sig_k^{-1}X^k -C) )}
\\[5pt]
&\leq&  \norm{\tZ^{k}-Z^{k+1} },
\end{eqnarray*}
where  $\tZ^k$ is the $Z$ at the penultimate iteration when we compute  $(y^{k+1},S^{k+1},Z^{k+1})=$ MSNCG $(y^k,S^k,Z^k,X^k,\sigma_k)$.
Note that
$\cA^* y^{k+1} + S^{k+1} +\tZ^{k} + \sig_k^{-1}X^k -C
= \Pi_{\cK}(\cA^*y^{k+1}+\tZ^k+\sig_k^{-1}X^k-C) \in \cK$.
Thus for a given $\zeta \in (0,+\infty)$, we can replace (B3) by the following criterion for terminating Algorithm MSNCG:
\begin{itemize}
\item[(B3$^{\prime}$)]
$\max\Big\{\norm{\cA (X^k+\sig_k R_D^{k+1}  )-b} , \zeta \sig_k \norm{Z^{k+1} -\tZ^{k}}\Big\} \leq (\delta^{'}_k/\sigma_k)  \|X^{k+1} -
X^k\|$,  $0 \leq \delta^{'}_k \downarrow 0$.
\end{itemize}

In our numerical experiments,  we measure the accuracy of an approximate optimal solution $(X,y,S,Z)$ for (P) and (D) by using the following relative residual:
\begin{eqnarray}
\eta = \max \{\eta_P,\eta_D,\eta_{\cK},\eta_{\cP},\eta_{\cK^*},\eta_{\cP^*},\eta_{C1},\eta_{C2}\},  \label{eq-eta}
\end{eqnarray}
where
$\eta_P =\frac{\norm{\cA X-b}}{1+\norm{b}}$,
$\eta_D =\frac{\norm{\cA^* y+S +Z-C}}{1+\norm{C}}$,
$\eta_{\cK} = \frac{\norm{\Pi_{\cK^*}(-X)}}{1+\norm{X}}$,
$\eta_{\cP} = \frac{\norm{\Pi_{\cP^*}(-X)}}{1+\norm{X}}$,
$\eta_{\cK^*} = \frac{\norm{\Pi_{\cK}(-S)}}{1+\norm{S}}$,
$\eta_{\cP^*} = \frac{\norm{\Pi_{\cP}(-Z)}}{1+\norm{Z}}$,
$\eta_{C1} = \frac{\abs{\inprod{X}{S}}}{1+\norm{X}+\norm{S}}$,
$\eta_{C2} = \frac{\abs{\inprod{X}{Z}}}{1+\norm{X}+\norm{Z}}$.
Additionally, we compute the relative gap by
\begin{eqnarray}
&\eta_g = \frac{\inprod{C}{X} - \inprod{b}{y} }{1+\abs{\inprod{C}{X}}+\abs{\inprod{b}{y}}}.&
\label{eq-gap}
\end{eqnarray}
Let  $\varepsilon > 0$ be a given  accuracy tolerance. We terminate both SDPNAL$+$ and  ADMM$+$   when
\begin{eqnarray}
\eta  < \varepsilon. \label{stop}
\end{eqnarray}

Note that SDPAD  can be used to solve  SDP$+$ problems of  form (P) with $\cP = \cS^n_{\geq 0}$ directly and we stop SDPAD when
$\eta < \varepsilon,$
where $\eta$ is defined as in (\ref{eq-eta}).
However, it is shown recently that the direct extension of ADMM to the  multi-block case is not necessarily convergent \cite{chen2013direct}. Hence SDPAD, which is essentially an implementation of the direct extension of ADMM with the step length set at $1.618$ for solving the dual of SDP$+$ problems,  does not have  convergence guarantee in theory.

The implementation of 2EBD including its termination, along with    ADMM$+$ and  SDPAD, is done in the same way as   in \cite{sun2014convergent}. For  2EBD, we reformulate QAP, RCP and R1TA problems as SDP problems in the standard form as these problems do not appear to have obvious two-easy blocks structures.

In our numerical experiments, we also use a restart strategy for SDPNAL$+$ if it is not able to achieve the required accuracy for   the tested SDP$+$ problems.  For some problems, even though $\eta_P$ and $\eta_D$ can reach the required accuracy tolerance, $\eta_{\cK}$ or $\eta_{C1}$ may stay above the required tolerance or stagnate. This may happen, as in the case for SDPNAL, because many of these SDP$+$ problems are degenerate at the optimal solutions. One way to overcome this difficulty is to apply ADMM$+$ to (P) using
the most recently computed $(y,S,Z,X,\sigma)$ to restart SDPNAL$+$ when its progress
is not satisfactory. From this point of view, our proposed  algorithm is quite flexible.

Table \ref{table-numsolved} shows the number  of problems that  have been successfully solved to the accuracy of  $10^{-6}$ in $\eta$ by each of  the four solvers SDPNAL$+$, ADMM$+$, SDPAD and 2EBD, with the maximum number of iterations set at $25000$ or the maximum
computation time set at $99$ hours. As can be seen, only SDPNAL$+$ can solve all the problems to the accuracy of  $10^{-6}$. In particular, for the first time, we are able to solve all the $95$ difficult SDP$+$ problems arising from QAP problems to an accuracy of $10^{-6}$ efficiently, while ADMM$+$, SDPAD and 2EBD can successfully solve 39, 30 and 16 problems, respectively.

\begin{table}[h]
\begin{center}
\caption{Number of problems which are solved to the accuracy of  $ 10^{-6}$ in $\eta$.} \label{table-numsolved}
\begin{tabular}{| c | c | c | c | c |}
\hline
\multicolumn{1}{|c}{problem set (No.) $\backslash$ solver} & \multicolumn{1}{|c|}{SDPNAL$+$} & \multicolumn{1}{|c|}{ADMM$+$} & \multicolumn{1}{|c|}{SDPAD}  & \multicolumn{1}{|c|}{2EBD}
 \\ \hline
$\theta$ (58) & 58 & 56 & 53 & 53 \\ \hline $\theta_{+}$ (58) & 58 & 58 & 58 & 56 \\ \hline FAP ( 7) &  7 &  7 &  7 &  7 \\ \hline QAP (95) & 95 & 39 & 30 & 16 \\ \hline BIQ (134) & 134 & 134 & 134 & 134 \\ \hline RCP (120) & 120 & 120 & 114 & 109 \\ \hline R1TA (55) & 55 & 42 & 47 & 18 \\ \hline Total (527) & 527 & 456 & 443 & 393 \\ \hline 
\hline
\end{tabular}
\end{center}
\end{table}

Tables \ref{table-SDPNALP-1} and \ref{table-SDPNALP-0} show the numerical results of SDPNAL$+$ with the tolerance $\varepsilon = 10^{-6}$. The first three columns of each table give the problem name, the dimension of the variable $y$ ($m$), the size of the matrix $C$ ($n_s$) and the number of linear inequality constraints ($n_l$) in (\PP), respectively. The middle five columns give the number of outer iterations, the total number of inner iterations, the total number of iterations for ADMM$+$, and the objective values $\inprod{C}{X}$ and $\inprod{b}{y}$. The relative infeasibilities and gap, as well as times (in the format hours:minutes:seconds) are listed in the last eight columns. It is interesting to note that all the tested problems (especially the QAPs) can be solved to the required accuracy $10^{-6}$ efficiently by SDPNAL$+$.

Tables \ref{tableLS_SDPNALP_ADMM3c_SDPAD_TEBD-1} and \ref{tableLS_SDPNALP_ADMM3c_SDPAD_TEBD-0} compare  SDPNAL$+$, ADMM$+$, SDPAD and 2EBD on a collection of SDP$+$ and SDP problems using the tolerance $\varepsilon = 10^{-6}$. We terminate ADMM$+$, SDPAD, 2EBD after 25000 iterations or 99 hours.
 As can be seen, except for SDPNAL$+$, the required accuracy is not achieved for most of the tested QAPs after 25000 iterations for the solvers ADMM$+$, SDPAD and 2EBD.
For the last three solvers, they typically converge very slowly
when $\eta$ falls below the range of $10^{-4}$--$10^{-5}$.
For R1TA problems, SDPNAL$+$ is significantly faster than the other 3 methods and it seems that only SDPNAL$+$ can solve those large scale ($n \geq 2000$) problems efficiently.

\begin{figure}[]
\begin{center}
\ifpdf
 \includegraphics[width=0.9\textwidth]{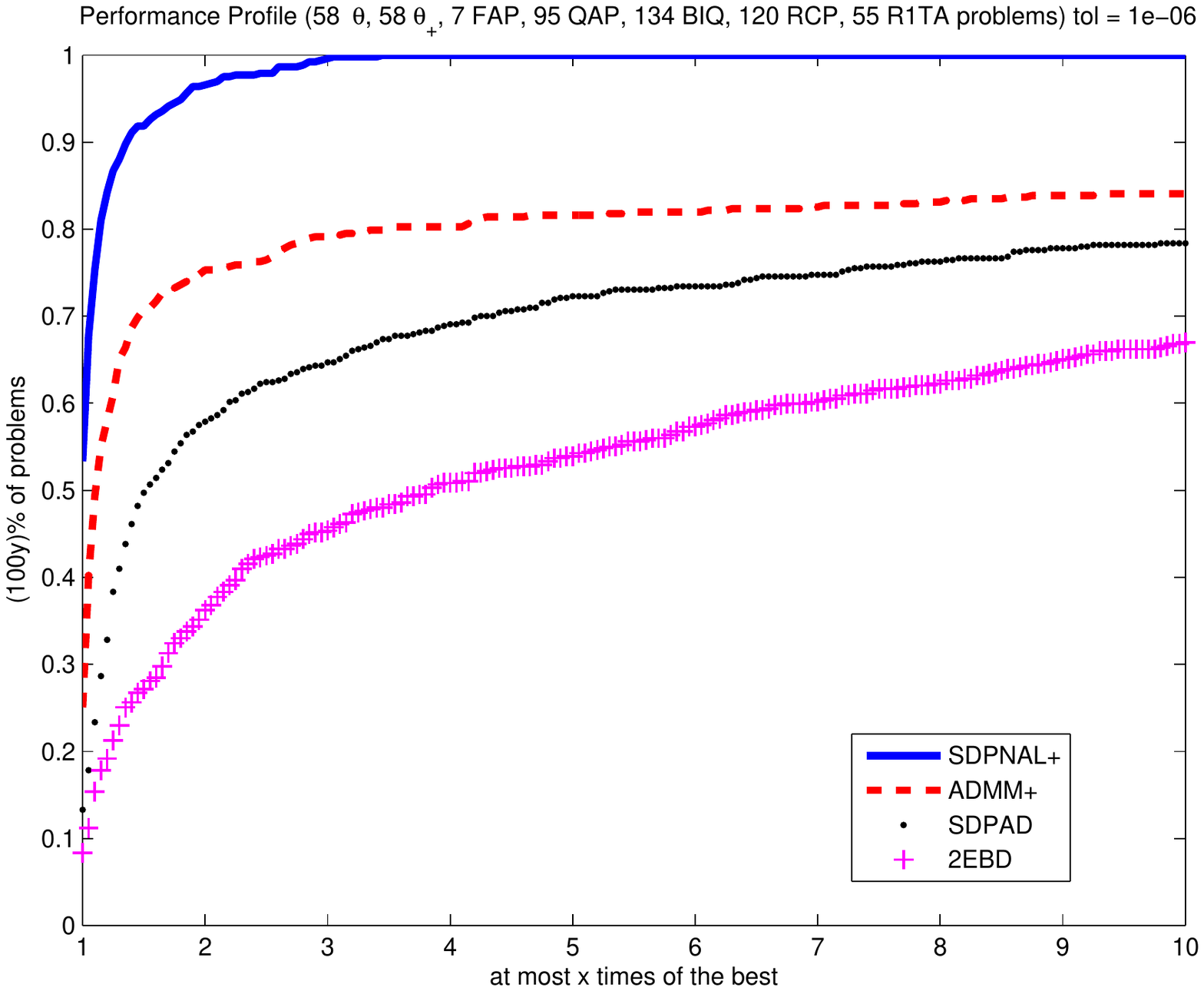}
 \else
 \includegraphics[width=0.9\textwidth]{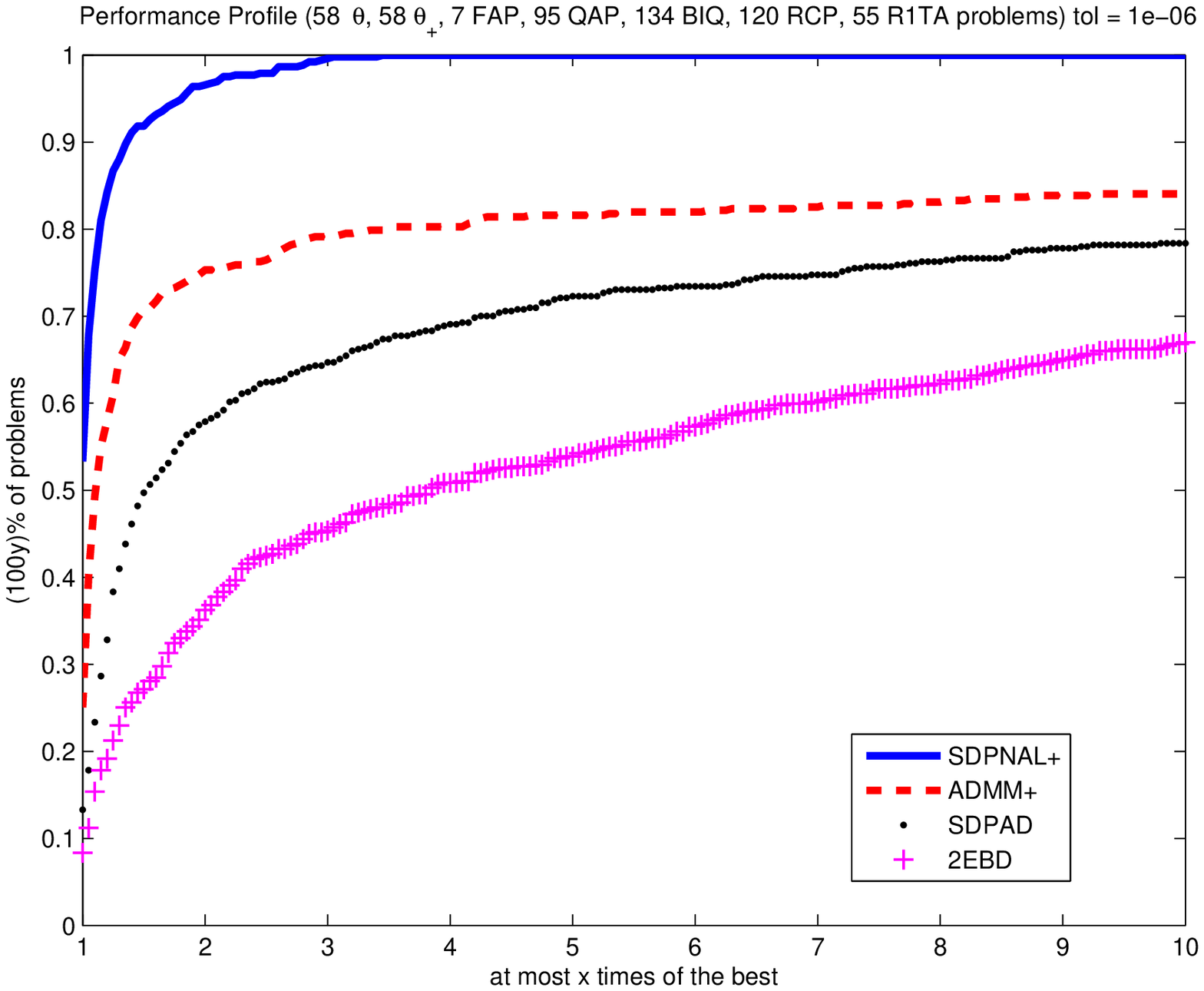}
 \fi
\caption{Performance profiles of SDPNAL$+$, ADMM$+$, SDPAD and 2EBD on $[1,10]$}
\label{fig-all}
\end{center}
\end{figure}

Figure \ref{fig-all} shows the performance profiles of SDPNAL$+$, ADMM$+$, SDPAD and 2EBD for the tested problems listed in Tables \ref{tableLS_SDPNALP_ADMM3c_SDPAD_TEBD-1} and \ref{tableLS_SDPNALP_ADMM3c_SDPAD_TEBD-0}. We recall that a point $(x,y)$ is in the performance profile curve of a method if and only if it can be solved exactly $(100y)\%$ of all the tested problems at most $x$ times slower than any other method. It can be seen that SDPNAL$+$ outperforms the other 3 methods by a significant margin.

\section*{Acknowledgements}

The authors would like to thank Jiawang Nie and Li Wang 
for sharing their codes on semidefinite relaxations of rank-1 tensor 
approximation problems.

\bibliographystyle{siam}
\bibliography{SDPNALPbib}

\newpage

\oddsidemargin -1.5cm
\evensidemargin 0pt
\topmargin 0cm
\textwidth 6.6in
\textheight 8.7in

\begin{center}
\begin{scriptsize}
\begin{longtable}{| ll | lr | lr|}

\caption{Performance of SDPNAL$+$ on $\theta_+$, FAP, QAP, BIQ and RCP problems ($\varepsilon = 10^{-6}$)} \label{table-SDPNALP-1} \\
\hline
\mc{2}{|c|}{} & \mc{2}{c|}{}  & \mc{2}{c|}{} \\[-8pt]
\mc{2}{|c|}{problem $|$ $m$ $|$ $n_s;n_l$}
&\mc{1}{l}{it$|$itsub$|$itA} &\mc{1}{c}{$pobj$ $\,|$ $dobj$}
&\mc{1}{l}{$\eta_P$ $\,|$ $\eta_D$ $\,|$ $\eta_{\cK_1}$ $\,|$ $\eta_{\cK_2}$ $\,|$ $\eta_{C1}$ $\,|$ $\eta_{C2}$ $\,|$ $\eta_g$}
&\mc{1}{c|}{time}
\\[2pt] \hline
\endhead
\hline
\endfoot
 &\mc{1}{c|}{} &\mc{2}{c|}{} &\mc{2}{c|}{}\\[-5pt] \hline 
\mc{2}{|c|}{theta4 $|$ 1949 $|$  200;}&  0$|$  0$|$304& 4.98689141 1 $|$ 4.98690560 1&   3.3-12$|$   9.6-7$|$   1.0-7$|$   5.8-8$|$   1.1-7$|$   1.1-8$|$   {\green{ -1.4-6}} &06\\  
 \mc{2}{|c|}{theta42 $|$ 5986 $|$  200;}&  0$|$  0$|$179& 2.37382082 1 $|$ 2.37382041 1&   3.7-14$|$   9.6-7$|$   7.8-8$|$   3.5-8$|$   3.4-9$|$   1.6-9$|$   8.3-8 &03\\  
 \mc{2}{|c|}{theta6 $|$ 4375 $|$  300;}&  0$|$  0$|$316& 6.29616118 1 $|$ 6.29618686 1&   5.9-12$|$   8.5-7$|$   1.0-10$|$   9.6-9$|$   1.6-7$|$   2.4-8$|$   {\green{ -2.0-6}} &13\\  
 \mc{2}{|c|}{theta62 $|$ 13390 $|$  300;}&  0$|$  0$|$178& 2.93779365 1 $|$ 2.93779433 1&   1.5-13$|$   9.5-7$|$   2.8-8$|$   1.2-8$|$   5.1-9$|$   2.2-9$|$   -1.1-7 &08\\  
 \mc{2}{|c|}{theta8 $|$ 7905 $|$  400;}&  0$|$  0$|$316& 7.34075146 1 $|$ 7.34078891 1&   9.4-14$|$   9.7-7$|$   1.6-10$|$   1.2-8$|$   1.7-7$|$   2.1-8$|$   {\green{ -2.5-6}} &24\\  
 \mc{2}{|c|}{theta82 $|$ 23872 $|$  400;}&  0$|$  0$|$157& 3.40643277 1 $|$ 3.40643435 1&   6.5-14$|$   9.7-7$|$   2.1-8$|$   9.8-9$|$   1.6-8$|$   4.0-9$|$   -2.3-7 &13\\  
 \mc{2}{|c|}{theta83 $|$ 39862 $|$  400;}&  0$|$  0$|$154& 2.01671025 1 $|$ 2.01671066 1&   4.7-14$|$   9.5-7$|$   2.8-8$|$   8.7-9$|$   3.3-9$|$   5.1-10$|$   -1.0-7 &14\\
\hline
\end{longtable}

\begin{longtable}{| ll | lr | lr|}
\caption{Performance of SDPNAL$+$ on $\theta$ and R1TA problems ($\varepsilon = 10^{-6}$)} \label{table-SDPNALP-0} \\
\hline
\mc{2}{|c|}{} & \mc{2}{c|}{}  & \mc{2}{c|}{} \\[-8pt]
\mc{2}{|c|}{problem $|$ $m$ $|$ $n_s;n_l$}
&\mc{1}{l}{it$|$itsub$|$itA} &\mc{1}{c}{$pobj$ $\,|$ $dobj$}
&\mc{1}{l}{$\eta_P$ $\,|$ $\eta_D$ $\,|$ $\eta_{\cK_1}$ $\,|$ $\eta_{\cK_2}$ $\,|$ $\eta_{C1}$ $\,|$ $\eta_{C2}$ $\,|$ $\eta_g$}
&\mc{1}{c|}{time}
\\[2pt] \hline
\endhead
\hline
\endfoot
 &\mc{1}{c|}{} &\mc{2}{c|}{} &\mc{2}{c|}{}\\[-5pt] \hline 
\mc{2}{|c|}{theta4 $|$ 1949 $|$  200;}& 13$|$ 13$|$153& 5.03212171 1 $|$ 5.03212139 1&   5.2-8$|$   6.5-7$|$   1.4-15$|$ 0 $|$   1.0-15$|$ 0 $|$   3.2-8 &04\\  
 \mc{2}{|c|}{theta42 $|$ 5986 $|$  200;}& 20$|$ 20$|$ 82& 2.39317014 1 $|$ 2.39317083 1&   2.5-7$|$   6.0-7$|$   1.1-15$|$ 0 $|$   6.7-16$|$ 0 $|$   -1.4-7 &05\\  
 \mc{2}{|c|}{theta6 $|$ 4375 $|$  300;}& 12$|$ 12$|$163& 6.34770795 1 $|$ 6.34770774 1&   3.7-8$|$   9.9-7$|$   3.4-16$|$ 0 $|$   1.2-15$|$ 0 $|$   1.6-8 &09\\  
 \mc{2}{|c|}{theta62 $|$ 13390 $|$  300;}& 13$|$ 13$|$ 82& 2.96412360 1 $|$ 2.96412502 1&   3.5-7$|$   6.7-7$|$   2.0-16$|$ 0 $|$   6.0-16$|$ 0 $|$   -2.4-7 &07\\  
 \mc{2}{|c|}{theta8 $|$ 7905 $|$  400;}& 12$|$ 12$|$183& 7.39535647 1 $|$ 7.39535733 1&   2.4-8$|$   4.4-7$|$   1.5-15$|$ 0 $|$   0.8-15$|$ 0 $|$   -5.8-8 &18\\  
 \mc{2}{|c|}{theta82 $|$ 23872 $|$  400;}& 11$|$ 11$|$ 87& 3.43668856 1 $|$ 3.43668909 1&   1.1-7$|$   4.5-7$|$   5.4-16$|$ 0 $|$   1.5-16$|$ 0 $|$   -7.6-8 &13\\  
 \mc{2}{|c|}{theta83 $|$ 39862 $|$  400;}& 23$|$ 23$|$ 64& 2.03018725 1 $|$ 2.03018894 1&   4.4-7$|$   9.8-7$|$   1.2-16$|$ 0 $|$   6.7-16$|$ 0 $|$   -4.1-7 &23\\
\hline
\end{longtable}
\end{scriptsize}
\end{center}

\begin{landscape}
\begin{scriptsize}
\begin{center}

\begin{longtable}{| c | c | c | c | c |}
\caption{Performance of SDPNAL$+$ (a), ADMM$+$ (b), SDPAD (c) and 2EBD (d) on $\theta_+$, FAP, QAP, BIQ and RCP problems ($\varepsilon = 10^{-6}$)} \label{tableLS_SDPNALP_ADMM3c_SDPAD_TEBD-1} \\
\hline
\multicolumn{1}{|c}{} & \multicolumn{1}{|c|}{iteration} & \multicolumn{1}{|c|}{$\eta$} & \multicolumn{1}{|c|}{$\eta_g$}  & \multicolumn{1}{|c|}{time} \\ \hline
problem $|$ $m$ $|$ $n_{s};n_{l}$
&a$|$ b$|$ c$|$ d
&a$|$ b$|$ c$|$ d
&a$|$ b$|$ c$|$ d
&a$|$ b$|$ c$|$ d
\\ \hline
\endhead
theta4 $|$  1949 $|$   200; &   0, 0;304 $|$  304 $|$  307 $|$  701& 9.6-7  $|$  9.6-7  $|$  9.9-7  $|$  9.9-7  &  {\green{ -1.4-6}}  $|$  {\green{ -1.4-6}}  $|$  3.0-7  $|$  {\green{ 1.5-6}}  & 06 $|$ 06 $|$ 05 $|$ 13  \\ \hline 
 theta42 $|$  5986 $|$   200; &   0, 0;179 $|$  179 $|$  151 $|$  345& 9.6-7  $|$  9.6-7  $|$  9.6-7  $|$  9.8-7  &  8.3-8  $|$  8.3-8  $|$  3.3-7  $|$  9.5-7  & 03 $|$ 03 $|$ 03 $|$ 06  \\ \hline 
 theta6 $|$  4375 $|$   300; &   0, 0;316 $|$  316 $|$  318 $|$  671& 8.5-7  $|$  8.5-7  $|$  9.8-7  $|$  9.9-7  &  {\green{ -2.0-6}}  $|$  {\green{ -2.0-6}}  $|$  {\green{ -1.2-6}}  $|$  {\green{ 1.4-6}}  & 13 $|$ 13 $|$ 12 $|$ 31  \\ \hline 
 theta62 $|$  13390 $|$   300; &   0, 0;178 $|$  178 $|$  129 $|$  336& 9.5-7  $|$  9.5-7  $|$  8.9-7  $|$  9.9-7  &  -1.1-7  $|$  -1.1-7  $|$  {\green{ 1.3-6}}  $|$  {\green{ 1.4-6}}  & 08 $|$ 08 $|$ 07 $|$ 15  \\ \hline 
 theta8 $|$  7905 $|$   400; &   0, 0;316 $|$  316 $|$  325 $|$  470& 9.7-7  $|$  9.7-7  $|$  8.6-7  $|$  9.8-7  &  {\green{ -2.5-6}}  $|$  {\green{ -2.5-6}}  $|$  {\green{ 1.7-6}}  $|$  {\green{ 1.2-6}}  & 24 $|$ 24 $|$ 24 $|$ 40  \\ \hline 
 theta82 $|$  23872 $|$   400; &   0, 0;157 $|$  157 $|$  130 $|$  343& 9.7-7  $|$  9.7-7  $|$  8.4-7  $|$  9.7-7  &  -2.3-7  $|$  -2.3-7  $|$  {\green{ 1.8-6}}  $|$  {\green{ 1.8-6}}  & 13 $|$ 13 $|$ 13 $|$ 30  \\ \hline 
 theta83 $|$  39862 $|$   400; &   0, 0;154 $|$  154 $|$  111 $|$  306& 9.5-7  $|$  9.5-7  $|$  9.9-7  $|$  9.9-7  &  -1.0-7  $|$  -1.0-7  $|$  6.0-7  $|$  {\green{ 2.4-6}}  & 14 $|$ 14 $|$ 13 $|$ 28  \\ \hline 
 theta10 $|$  12470 $|$   500; &   0, 0;354 $|$  354 $|$  351 $|$  490& 8.5-7  $|$  8.5-7  $|$  9.9-7  $|$  9.9-7  &  {\green{ -2.5-6}}  $|$  {\green{ -2.5-6}}  $|$  -7.8-7  $|$  {\green{ 1.5-6}}  & 46 $|$ 45 $|$ 42 $|$ 1:11  \\
\hline 
 
\end{longtable}

\begin{longtable}{| c | c | c | c | c |}
\caption{Performance of SDPNAL$+$ (a), ADMM$+$ (b), SDPAD (c) and 2EBD (d)
on  $\theta$ and R1TA problems ($\varepsilon = 10^{-6}$)} \label{tableLS_SDPNALP_ADMM3c_SDPAD_TEBD-0} \\
\hline
\multicolumn{1}{|c}{} & \multicolumn{1}{|c|}{iteration} & \multicolumn{1}{|c|}{$\eta$} & \multicolumn{1}{|c|}{$\eta_g$}  & \multicolumn{1}{|c|}{time} \\ \hline
problem $|$ $m$ $|$ $n_{s};n_{l}$
&a$|$ b$|$ c$|$ d
&a$|$ b$|$ c$|$ d
&a$|$ b$|$ c$|$ d
&a$|$ b$|$ c$|$ d
\\ \hline
\endhead
theta4 $|$  1949 $|$   200; &  13,13;153 $|$  316 $|$  294 $|$  680& 6.5-7  $|$  9.7-7  $|$  9.6-7  $|$  9.9-7  &  3.2-8  $|$  {\green{ -1.3-6}}  $|$  -5.9-7  $|$  {\green{ 1.5-6}}  & 04 $|$ 05 $|$ 04 $|$ 12  \\ \hline 
 theta42 $|$  5986 $|$   200; &  20,20;82 $|$  178 $|$  151 $|$  315& 6.0-7  $|$  9.9-7  $|$  9.9-7  $|$  9.8-7  &  -1.4-7  $|$  -1.0-7  $|$  2.7-7  $|$  6.3-7  & 05 $|$ 03 $|$ 03 $|$ 06  \\ \hline 
 theta6 $|$  4375 $|$   300; &  12,12;163 $|$  341 $|$  306 $|$  561& 9.9-7  $|$  8.4-7  $|$  9.5-7  $|$  9.9-7  &  1.6-8  $|$  {\green{ 1.7-6}}  $|$  {\green{ -1.1-6}}  $|$  {\green{ 1.7-6}}  & 09 $|$ 13 $|$ 10 $|$ 24  \\ \hline 
 theta62 $|$  13390 $|$   300; &  13,13;82 $|$  190 $|$  125 $|$  311& 6.7-7  $|$  9.7-7  $|$  9.6-7  $|$  9.9-7  &  -2.4-7  $|$  7.5-9  $|$  {\green{ 1.1-6}}  $|$  {\green{ 1.2-6}}  & 07 $|$ 08 $|$ 06 $|$ 14  \\ \hline 
 theta8 $|$  7905 $|$   400; &  12,12;183 $|$  364 $|$  321 $|$  413& 4.4-7  $|$  9.9-7  $|$  8.3-7  $|$  9.9-7  &  -5.8-8  $|$  {\green{ -2.1-6}}  $|$  9.8-7  $|$  {\green{ 1.0-6}}  & 18 $|$ 28 $|$ 21 $|$ 34  \\ \hline 
 theta82 $|$  23872 $|$   400; &  11,11;87 $|$  163 $|$  126 $|$  304& 4.5-7  $|$  9.8-7  $|$  9.7-7  $|$  9.7-7  &  -7.6-8  $|$  3.8-7  $|$  {\green{ 2.0-6}}  $|$  {\green{ 1.3-6}}  & 13 $|$ 14 $|$ 12 $|$ 26  \\ \hline 
 theta83 $|$  39862 $|$   400; &  23,23;64 $|$  157 $|$  110 $|$  290& 9.8-7  $|$  9.2-7  $|$  9.2-7  $|$  9.8-7  &  -4.1-7  $|$  -1.1-7  $|$  2.7-7  $|$  {\green{ 1.8-6}}  & 23 $|$ 14 $|$ 13 $|$ 26  \\ \hline 
 theta10 $|$  12470 $|$   500; &  11,11;200 $|$  396 $|$  333 $|$  422& 7.6-7  $|$  9.1-7  $|$  9.9-7  $|$  9.8-7  &  6.7-8  $|$  {\green{ -2.1-6}}  $|$  {\green{ -1.4-6}}  $|$  {\green{ 1.2-6}}  & 32 $|$ 51 $|$ 36 $|$ 59  \\
\hline 
 
\end{longtable}
\end{center}
\end{scriptsize}
\end{landscape}

\end{document}